\documentclass[a4paper,11pt]{amsart}
\usepackage[letterpaper, margin=1.2in]{geometry}
\usepackage{latexsym}
\usepackage{amssymb}
\usepackage{amsmath}
\usepackage{amsbsy}
\usepackage{amsfonts}
\usepackage[table]{xcolor}
\usepackage{pgfplots}
\usepackage{color}
\usepackage{multirow}
\usepackage{graphicx}
\usepackage{txfonts,pxfonts} 
\usepackage{tikz}
\usetikzlibrary{calc,arrows}
\usetikzlibrary{calc}
\usepackage{verbatim}
\everymath{\displaystyle}
\newtheorem{thm}{Theorem}[section]

\newtheorem{lem}[thm]{Lemma}

\newtheorem{cor}[thm]{Corollary}

\theoremstyle{definition}

\newcommand{\Z}{\mathbb Z}

\def\ol{\overline}

\def\la{\lambda}

\def\md#1{\ \mbox{\rm(mod }{#1})}

\newcounter{cs}
\stepcounter{cs}
\newcommand{\casos}{\begin{itemize}}
	\newcommand{\fcasos}{\end{itemize}\setcounter{cs}{1}}

\newfont{\tit}{cmr12 scaled \magstep3}
\title[The index of certain nonic number fields defined by $x^9+ax^2+b$]{The index of certain nonic number fields defined by $x^9+ax^2+b$}
\author[L. EL Fadil]{Lhoussain El Fadil}
\author[O. Kchit]{Omar Kchit}
\address{Faculty of Sciences Dhar El Mahraz, P.O. Box  1796 Atlas-Fes, Sidi Mohamed ben Abdellah University,  Morocco}
\email{lhoussain.elfadil@usmba.ac.ma, \,\, orcid: 0000.0003.4175.8064} 
\email{omar.kchit@usmba.ac.ma, \,\, orcid: 0000.0002.0844.5034}
\begin{document}
	
	\begin{abstract}
		In this	paper, for any nonic number field $K$ defined by a monic irreducible trinomial $F(x)=x^9+ax^2+b \in \mathbb{Z}[x]$, we calculate $\nu_p(i(K))$ for every rational prime $p$. In particular, we characterize the index $i(K)$ of this family of number fields. As an application of our results, if $i(K)\neq1$, then $K$ is not monogenic. We illustrate our results by some computational examples. 
	\end{abstract}
	\maketitle
	\section{Introduction}
	Let $K$ be a number field of degree $n$ with ring of integers $\mathbb{Z}_K$. Let $\alpha\in \mathbb{Z}_K$ be a primitive element of $K$. The index of $\alpha$, denoted by $(\mathbb{Z}_K: \mathbb{Z}[\alpha])$, is the index the Abelian group $\mathbb{Z}[\alpha]$ in $\mathbb{Z}_K$, and the formula linking $(\mathbb{Z}_K: \mathbb{Z}[\alpha])$, $\Delta(\alpha)$, and $d_K$ is given by:
	\begin{equation}\tag{1.1}
		\Delta(\alpha)=\pm (\mathbb{Z}_K: \mathbb{Z}[\alpha])^2\cdot d_K,
	\end{equation}  
	where $d_K$ is the absolute discriminant of $K$ and $\Delta(\alpha)$ is the discriminant of the minimal polynomial of $\alpha$ over $\mathbb{Q}$. The index of $K$, denoted by $i(K)$, is the greatest common divisor of the indices of all primitive integers of $K$. Say $i(K)=\gcd \ \{ ( \mathbb{Z}_K:\mathbb{Z}[\theta]) \, |\,{K=\mathbb{Q}(\theta)\text{ and }\theta \in \mathbb{Z}_K}\}$. A rational prime $p$ dividing $i(K)$ is called a prime common index divisor of $K$. If $K$ is monogenic, then $\mathbb{Z}_K$ has a power integral basis, i.e., a $\mathbb{Z}$-basis of the form $(1,\theta,\dots,\theta^{n-1})$, and the index of $K$ is trivial, say $i(K)=1$. Therefore a field having a prime common index divisor is not monogenic. {The first number field with non trivial index was given by Dedekind in $1871$, who exhibited the example of the cubic number field $K$ defined by $x^3-x^2-2x-8$ and showed that the index of any primitive integer of $K$ is even. In $1930$, Engstrom} \cite{En} was the first one who studied the prime power decomposition of the index of a number field. He showed that, for number fields of degree $n\leq 7$,  $\nu_p(i(K))$ is determined by the form of the factorization of $p\mathbb{Z}_K$, where $\nu_p$ is the $p$-adic valuation of $\mathbb{Q}$. This motivated a very important question, stated as problem $22$ in Narkiewicz's book (\cite{Nar}), which asks for an explicit formula of the highest power $\nu_p(i(K))$ for a given rational prime p dividing $i(K)$. In \cite{Sl}, \'Sliwa showed that, if $p$ is a non-ramified ideal in $K$, then $\nu_p(i(K))$ is determined by the factorization of $p\mathbb{Z}_K$. These results were generalized by Nart (\cite{N}), who developed a $p$-adic characterization of the index of a number field. In \cite{Nak}, Nakahara studied the index of non-cyclic but abelian biquadratic number fields. In \cite{GPP}  Ga\'al et al. characterized the field indices of biquadratic number fields having Galois group $V_4$. In \cite{EK6}, we characterized the index and studied the monogenity of any sextic number field defined by $x^6+ax^5+b$. In \cite{EK7}, for every rational prime $p$, we characterized $\nu_p(i(K))$ for any septic number field defined by a trinomial $x^7+ax^3+b$. In \cite{E53}, El Fadil characterized the index of any quintic number field defined by $x^5+ax^3+b$. In \cite{K9}, Kchit studied the index of any nonic number field defined by $x^9+ax+b$. In this paper, for any nonic number field $K$ defined by a monic irreducible trinomial $x^9+ax^2+b\in \mathbb{Z}[x]$, we calculate $\nu_p(i(K))$ for every rational prime $p$. In particular, we show that $i(K)\in\{0,1,2,3,6\}$.
	\section{Main Results}
	Throughout this section, $K$ is a number field generated by a complex root $\alpha$ of a {monic irreducible} trinomial $F(x)=x^9+ax^2+b\in\mathbb{Z}[x]$. Without loss of generality, we assume that for every rational prime $p$,  $\nu_p(a)\leq 6$ or  $\nu_p(b) \leq 8$.\\
	Along this paper, for every integer $m\in\mathbb{Z}$ and a rational prime $p$, let $m_p=\cfrac{m}{p^{\nu_p(m)}}$. We denote $\Delta=\Delta(F)$ the discriminant of $F(x)$.\\
	\smallskip
	The following theorem characterizes when $2$ divides $i(K)$. 
	\begin{thm}\label{thmp2}
		The following table provides $\nu_2(i(K))$:
		
		\begin{table}[htbp]
			\centering
			\caption{$\nu_2(i(K))$}
			\begin{tabular}{|l|c|c|}
				\hline
				\multicolumn{2}{|c|}{\textbf{Conditions}}&\textbf{$\quad\nu_2(i(K))\quad$}\\
				\hline
				$\nu_2(a)=0$ and $\nu_2(b)=2k$ &$b_2\equiv -a\md{8}$&\multirow{9}{*}{$1$}\\
				\cline{1-2}  
				$a\equiv 2\md{4}$ and $\nu_2(b)=2k+1$&\multirow{6}{*}{$b_2\equiv -a_2\md{8}$}&\\
				\cline{1-1}
				$a\equiv 4\md{8}$ and $\nu_2(b)=2k+2$&&\\
				\cline{1-1}
				$a\equiv 8\md{16}$ and $\nu_2(b)=2k+3$&&\\
				\cline{1-1}
				$a\equiv 16\md{32}$ and $\nu_2(b)=2k+4$&&\\
				\cline{1-1}
				$a\equiv 32\md{64}$ and $\nu_2(b)=2k+7$&&\\
				\cline{1-1}
				$a\equiv 64\md{128}$ and $\nu_2(b)=2k+8$&&\\
				\cline{1-2}
				$a\equiv 32\md{64}$ and $b\equiv 128\md{256}$&$a_2+b_2\equiv 4\md{8}$&\\
				\cline{1-2}
				$a\equiv 64\md{128}$ and $b\equiv 256\md{512}$&$\nu_2(\Delta)=2k$ and $\Delta_2\equiv a_2b_2\md{8}$&\\
				\hline
				\multicolumn{2}{|c|}{\textbf{Otherwise}}&$0$\\
				\hline
			\end{tabular}
		\end{table}
	\end{thm}
	
	\smallskip
	The following theorem characterizes when $3$ divides $i(K)$.
	\begin{thm}\label{thmp3} 
		The following table provides $\nu_3(i(K))$:
		
	\begin{table}[htbp]
		\centering
		\caption{$\nu_3(i(K))$}
		\begin{tabular}{|l|c|c|}
			\hline
			\multicolumn{2}{|c|}{\textbf{Conditions}}&\textbf{$\nu_3(i(K))$}\\
			\hline
			$(a,b)\in\{(9,71),(36,44)\}\md{81}$ and $a+b\equiv -1\md{243}$&$\nu_3(\Delta)=2k$&\multirow{4}{*}{$1$}\\
			\cline{1-1}
			$(a,b)\equiv(63,17)\md{81}$ and $a+b\equiv 80\md{243}$&\multirow{2}{*}{and}&\\
			\cline{1-1} 
			$(a,b)\equiv(18,64)\md{81}$ and $a+b\equiv 163\md{243}$&&\\
			\cline{1-1} 
			$(a,b)\in\{(45,37),(72,10)\}\md{81}$ and  $a+b\equiv 1\md{243}$& $\Delta_3\equiv -1\md{3}$&\\
			\hline
			\multicolumn{2}{|c|}{\textbf{Otherwise}}&$0$\\
			\hline
		\end{tabular}
	\end{table}
	\end{thm}
	\smallskip 
	
	\begin{thm}\label{thmp}
		For every rational prime $p\geq 5$ and {for} every $(a,b)\in\mathbb{Z}^2$ such that $F(x)=x^9+ax^2+b$ is irreducible, $p$ does not divide the index $i(K)$, where $K$ is the number field defined by $F(x)$.
	\end{thm}
	\smallskip
	\begin{cor}
		For every integers $a$ and $b$ such that $F(x)=x^9+ax^2+b$ is irreducible over $\mathbb{Q}$, $i(K)\in\{1,2,3,6\}$, where $K$ is the number field generated by a root of $F(x)$.
	\end{cor}

	\section{Preliminaries}\label{sec3}
	{Our proofs are based on Newton polygon techniques applied on prime ideal factorization, which is a standard method which is rather technical but very efficient to apply. We have introduced the corresponding concepts in several former papers. Here we only give the theorem of index of Ore which plays a key role for proving our main results.}\\
	Let $K=\mathbb{Q}(\alpha)$ be a number field generated by a complex root $\alpha$ of a monic irreducible polynomial $F(x)\in\mathbb{Z}[x]$. We shall use Dedekind’s theorem \cite[Chapter I, Proposition 8.3]{Neu} and Dedekind’s criterion \cite[Theorem 6.1.4]{Co}. Let $\phi\in\mathbb{Z}_p[x]$ be a monic lift to an irreducible factor of $F(x)$ modulo $p$, $F(x)=a_0(x)+a_1(x)\phi(x)+\cdots+a_l(x)\phi(x)^l$ the $\phi$-expansion of $F(x)$, and $N_{\phi}^+(F)$ the principal $\phi$-Newton polygon of $F(x)$. Let $\mathbb{F}_{\phi}$ be the field $\mathbb{F}_p[x]/(\overline{\phi})$, then to every side $S$ of $N_{\phi}^+(F)$ with initial point $(i,u_i)$, and every $i=0,\ldots,l$, let the residue coefficient $c_i\in\mathbb{F}_{\phi}$ defined as follows: 
	$$c_{i}=
	\left
	\{\begin{array}{ll} 0,& \mbox{ if } (s+i,{\it u_{s+i}}) \mbox{ lies strictly
			above } S,\\
		\left(\dfrac{a_{s+i}(x)}{p^{{\it u_{s+i}}}}\right)
		\,\,
		\mod{(p,\phi(x))},&\mbox{ if }(s+i,{\it u_{s+i}}) \mbox{ lies on }S.
	\end{array}
	\right.$$
	Let $-\lambda=-h/e$ be the slope of $S$, where $h$ and $e$ are two positive coprime integers and $l=l(S)$ its length. Then  $d=l/e$ is the degree of $S$. Hence, if $i$ is not a multiple of $e$, then  $(s+i, u_{s+i})$ does not lie on $S$, and so $c_i=0$. Let ${R_{\lambda}(F)(y)}=t_dy^d+t_{d-1}y^{d-1}+\cdots+t_{1}y+t_{0}\in\mathbb{F}_{\phi}[y]$, called  
	the residual polynomial of $F(x)$ associated to the side $S$, where for every $i=0,\dots,d$,  $t_i=c_{ie}$. If ${R_{\lambda}(F)(y)}$ is square free for each side of the polygon $N_{\phi}^+(F)$, then we say that $F(x)$ is $\phi$-regular.\\ 
	Let $\overline{F(x)}=\prod_{i=1}^{r}\ol{\phi_i}^{l_i}$ be the factorization of $F(x)$ into powers of monic irreducible coprime polynomials over $\mathbb{F}_p$, we say that the polynomial $F(x)$ is $p$-{regular} if $F(x)$ is a $\phi_i$-regular polynomial with respect to $p$ for every $i=1,\dots,r$. Let  $N_{\phi_i}^+(F)=S_{i1}+\cdots+S_{ir_i}$ be the $\phi_i$-principal Newton polygon of $F(x)$ with respect to $p$. For every $j=1,\dots,r_i$, let ${R_{\lambda_{ij}}(F)(y)}=\prod_{s=1}^{s_{ij}}\psi_{ijs}^{a_{ijs}}(y)$ be the factorization of {$R_{\lambda_{ij}}(F)(y)$} in $\mathbb{F}_{\phi_i}[y]$. Then we have the following  theorem of index of Ore:
	\begin{thm}\label{thm4}$($\cite[Theorem 1.7 and Theorem 1.9]{EMN}$)$\\
		Under the above hypothesis, we have the following:
		\begin{enumerate}
			\item 
			$$\nu_p((\mathbb{Z}_K:\mathbb{Z}[\alpha]))\geq\sum_{i=1}^{r}\text{ind}_{\phi_i}(F).$$  
			The equality holds if $F(x) \text{ is }p$-regular.
			\item 
			If $F(x) \text{ is }p$-regular, then
			$$p\mathbb{Z}_K=\prod_{i=1}^r\prod_{j=1}^{t_i}\prod_{s=1}^{s_{ij}}\mathfrak{p}_{ijs}^{e_{ij}}$$
			is the factorization of $p\mathbb{Z}_K$ into powers of prime ideals of $\mathbb{Z}_K$, where $e_{ij}$ is the smallest positive integer satisfying $e_{ij}\la_{ij}\in \Z$ and the residue degree of $\mathfrak{p}_{ijs}$ over $p$ is given by $f_{ijs}=\deg(\phi_i)\cdot \deg(\psi_{ijs})$ for every $(i,j,s)$.
		\end{enumerate}
	\end{thm} 
    \smallskip
    If the theorem of Ore fails, that is, $F(x)$ is not $p$-regular, then in order to complete the factorization of $F(x)$, Guàrdia, Montes, and Nart introduced the notion of high order Newton polygon (\cite{GMN}). Similar to first order, for each order $r$, they introduced the valuation $\omega_r$, the key polynomial $\phi_r(x)$ for this valuation, the Newton polygon $N_r(F)$ of any polynomial $F(x)$ with respect to $\omega_r$ and $\phi_r(x)$, and for each side $T_i$ of $N_r(F)$, the residual polynomial $R_{r}(F)(y)$, and the index of $F(x)$ in order $r$. For more details, we refer to \cite{GMN}.\\
	
	For the proof of our results, we need the following lemma, which characterizes the prime common index divisors of $K$.
	\begin{lem}\label{index}$($\cite{En}$)$\\
		Let $p$ be a rational prime and {$K$ a} number field. For every positive integer $f$, let $\mathcal{P}_f$ be the number of distinct prime ideals of $\mathbb{Z}_K$ lying above $p$ with residue degree $f$ and $\mathcal{N}_f$ the number of monic irreducible polynomials of $\mathbb{F}_p[x]$ of degree $f$.  Then $p$ {divides the index $i(K)$} if and only if $\mathcal{P}_f>\mathcal{N}_f$ for some positive integer $f$.
	\end{lem}
	
	\section{Proofs of Main Results}
	Throughout this section, if $p\mathbb{Z}_K=\prod_{i=1}^r\prod_{j=1}^{t_i}\prod_{s=1}^{s_{ij}}\mathfrak{p}_{ijs}^{e_{ij}}$, then $e_{ij}$ denotes the ramification index of $\mathfrak{p}_{ijs}$ and $f_{ijs}$ denotes its residue degree for every $(i,j,s)$.\\
	\smallskip
	
	\textit{\textbf{Proof of Theorem \ref{thmp2}}}.\\
	If $\nu_2(b)=0$, then since $\Delta=b(3^{18}b^7+2^{2}\cdot7^7a^9)$ is the discriminant of $F(x)$, thanks to the index formula $(1.1)$, $\nu_2((\mathbb{Z}_K:\mathbb{Z}[\alpha]))=0$ and so $\nu_2(i(K))=0$. Assume that $2$ divides $b$. Then we have the following cases:
	\begin{enumerate}
		\item If $\nu_2(a)=0$, then $F(x)\equiv x^2(x-1)(x^3+x+1)(x^3+x^2+1)\md{2}$. Then $2$ divides the index $i(K)$ if and only if $\phi=x$ provides two prime ideals of $\mathbb{Z}_K$ lying above $2$ with residue degree $1$ each ideal factor. That is $\nu_2(b)=2k$ and $b_2+a\equiv 0\md{8}$. In this case, we replace $\phi$ by $g=x-2^k$. Then $F(x)=\cdots+(a+36(2^{7k}))g^2+(2^{k+1}a+9(2^{8k}))g+b+2^{2k}a+2^{9k}$. Thus $N_g^+(F)=S_1+S_2$ has two sides joining $(0,w)$, $(1,k+1)$, and $(2,0)$ with $w\geq 2k+3$. Hence the degree of each side is $1$, and so $2\mathbb{Z}_K=\mathfrak{p}_{11}\mathfrak{p}_{12}\mathfrak{p}_{21}\mathfrak{p}_{31}\mathfrak{p}_{41}$ with $f_{11}=f_{12}=f_2=1$ and $f_3=f_4=3$. Applying \cite[Theorem 4]{En}, we get $\nu_2(i(K))=1$.

		\item If $\nu_2(a)\geq 1$, then $F(x)=x^9\md{2}$. Let $\phi=x$. By assumption, $\nu_2(a)\leq6$ or $\nu_2(b)\leq8$.\\
		If $\nu_2(b)<\cfrac{9}{7}\nu_2(a)$, then $N_{\phi}(F)=S_1$ has a single side joining $(0,\nu_2(b))$ and $(9,0)$ with $d(S_1)\in\{1,3\}$ since $\nu_2(b)\leq 8$ by assumption. Let $d=d(S_1)$, then we have the following cases:
		\begin{enumerate}
			\item[(i)] If $d=1$; $\nu_2(b)\in\{1,2,4,5,7,8\}$, then $2\mathbb{Z}_K=\mathfrak{p}_1^9$ with residue degree $1$. Hence $\nu_2(i(K))=0$.
			
			\item[(ii)] If $d=3$; $\nu_2(b)\in\{3,6\}$, then $R_{\lambda}(F)(y)=y^3+1=(y+1)(y^2+y+1)\in\mathbb{F}_{\phi}[y]$. Thus $2\mathbb{Z}_K=\mathfrak{p}_{11}^3\mathfrak{p}_{12}^3$ with $f_{11}=1$ and $f_{12}=2$. Hence $\nu_2(i(K))=0$.
		\end{enumerate}
		If $\nu_2(b)>\cfrac{9}{7}\nu_2(a)$, then $N_{\phi}(F)=S_1+S_2$ has two sides joining $(0,\nu_2(b))$, $(2,\nu_2(a))$, and $(9,0)$ with $d(S_1)\in\{1,2\}$ and $d(S_2)=1$ since $\nu_2(a)\leq 6$ by assumption. 
		\begin{enumerate}
			\item[(i)] If $\nu_2(b)\not\equiv \nu_2(a)\md{2}$, then $d(S_1)=1$. Thus $2\mathbb{Z}_K=\mathfrak{p}_1^2\mathfrak{p}_2^7$ with residue degree $1$ each ideal factor. Hence $\nu_2(i(K))=0$.
			
			\item[(ii)] If $\nu_2(b)=2k+\nu_2(a)$, then $d(S_1)=2$ with $R_{\lambda_1}(F)(y)=(y+1)^2\in\mathbb{F}_{\phi}[y]$. Let us replace $\phi$ by $g=x-2^k$. Then $F(x)=\cdots+(a+36(2^{7k}))g^2+(2^{k+1}a+9(2^{8k}))g+b+2^{2k}a+2^{9k}$. Then $\nu_2(a+36(2^{7k}))=\nu_2(a)$.\\
			If $k\geq2$ or $k=1$ and $\nu_2(a)\in\{1,2,3,4\}$, then $\nu_2(2^{k+1}a+9(2^{8k}))=\nu_2(a)+k+1$, and $\nu_2(b+2^{2k}a+2^{9k})\geq \nu_2(a)+2k+1$.
			\begin{enumerate}
				\item[(a)] If $b_2+a_2\equiv 2\md{4}$, then $N_g(F)=S_1+S_2$ has two sides joining $(0,2k+1+\nu_2(a))$, $(2,\nu_2(a))$, and $(9,0)$. Thus the degree of each side is $1$, and so $2\mathbb{Z}_K=\mathfrak{p}_1^2\mathfrak{p}_2^7$ with residue degree $1$ each ideal factor. Hence $\nu_2(i(K))=0$.
				
				\item[(b)] If $b_2+a_2\equiv 4\md{8}$, then $N_g(F)=S_1+S_2$ has two sides joining $(0,2k+2+\nu_2(a))$, $(2,\nu_2(a))$, and $(9,0)$ with $d(S_2)=1$ and $R_{\lambda_1}(F)=y^2+y+1$, which is irreducible over $\mathbb{F}_g$. Thus $2\mathbb{Z}_K=\mathfrak{p}_1\mathfrak{p}_2^7$ with $f_1=2$ and $f_2=1$. Hence $\nu_2(i(K))=0$.
				
				\item[(c)] If $b_2+a_2\equiv 0\md{8}$, then $N_g(F)=S_1+S_2+S_3$ has three sides joining $(0,w)$, $(1,k+1+\nu_2(a))$, $(2,\nu_2(a))$, and $(9,0)$ with $w\geq 2k+3+\nu_2(a)$. Thus the degree of each side is $1$, and so $2\mathbb{Z}_K=\mathfrak{p}_1\mathfrak{p}_2\mathfrak{p}_3^7$ with residue degree $1$ each ideal factor. Hence $2$ divides $i(K)$. Applying \cite[Corollary, p: 230]{En}, we get $\nu_2(i(K))=1$.
			\end{enumerate}
			If $k=1$ and $\nu_2(a)=5$; $a\equiv 32\md{64}$ and $b\equiv 128\md{256}$, then $\nu_2(4a+9(2^{8})=7$, and $\nu_2(b+4a+2^{9})\geq 8$.
	\begin{enumerate}
		\item[(a)] If $b_2+a_2\equiv 2\md{4}$, then $N_g(F)=S_1+S_2$ has two sides joining $(0,8)$, $(2,5)$, and $(9,0)$. Thus the degree of each side is $1$, and so $2\mathbb{Z}_K=\mathfrak{p}_1^2\mathfrak{p}_2^7$ with residue degree $1$ each ideal factor. Hence $\nu_2(i(K))=0$.
		
		\item[(b)] If $b_2+a_2\equiv 4\md{8}$, then $N_g(F)=S_1+S_2+S_3$ has three sides joining $(0,w)$, $(1,7)$, $(2,5)$, and $(9,0)$ with $w\geq 10$. Thus the degree of each side is $1$, and so $2\mathbb{Z}_K=\mathfrak{p}_1\mathfrak{p}_2\mathfrak{p}_3^7$ with residue degree $1$ each ideal factor. Hence $2$ divides $i(K)$. Applying \cite[Corollary, p: 230]{En}, we get $\nu_2(i(K))=1$.
		
		\item[(c)] If $b_2+a_2\equiv 0\md{8}$, then $N_g(F)=S_1+S_2$ has two sides joining $(0,9)$, $(2,5)$, and $(9,0)$ with $d(S_2)=1$ and $R_{\lambda_1}(F)=y^2+y+1$, which is irreducible over $\mathbb{F}_g$. Thus $2\mathbb{Z}_K=\mathfrak{p}_1\mathfrak{p}_2^7$ with $f_1=2$ and $f_2=1$. Hence $\nu_2(i(K))=0$.
		\end{enumerate}
		If $k=1$ and $\nu_2(a)=6$; $a\equiv 64\md{128}$ and $b\equiv 256\md{512}$, then $\nu_2(\Delta)\geq65$. Let $u=\cfrac{2\cdot 7^3a_2^4}{3^8b_2^3}$. Since $\nu_2=(3^8b_2^3)=0$, then $u\in\mathbb{Z}_2$. Let $\phi_1=x-u$. Then $F(x)=\cdots+(36u^7+a)\phi_1^2+A\phi_1+B$, where $\nu_2(36u^7+a)=6$,
		$$\begin{array}{lllll}
			A&=&9u^8+2au&=&\cfrac{3^7a_2^5\Delta^3}{2^{184}\cdot3^{62}b_2^{27}}-\cfrac{3^7a_2^5\Delta^2}{2^{120}\cdot3^{43}b_2^{19}}+\cfrac{3^7a_2^5\Delta}{2^{56}\cdot3^{25}b_2^{11}}~~,\text{ and}\\
			&&&&\\
			B&=&u^9+au^2+b&=&\cfrac{\Delta^4}{2^{246}\cdot3^{72}\cdot 7b_2^{31}}-\cfrac{\Delta^3}{2^{181}\cdot3^{54}\cdot 7b_2^{23}}+\cfrac{\Delta^2}{2^{118}\cdot3^{35}\cdot 7b_2^{15}}+\cfrac{\Delta}{2^{56}\cdot3^{18}\cdot7b_2^{7}}~~.
		\end{array}$$
	Thus $\nu_2(A)=\nu_2(B)=\nu_2(\Delta)-56$. It follows that $N_{\phi_1}(F)=S_1+S_2$ has two side joining $(0,\nu_2(\Delta)-56)$, $(2,6)$, and $(9,0)$ with $d(S_2)=1$. If $\nu_2(\Delta)$ is odd, then $d(S_1)=1$, and so $2\mathbb{Z}_K=\mathfrak{p}_1^2\mathfrak{p}_2^7$ with residue degree $1$ each ideal factor. Hence $\nu_2(i(K))=0$. If $\nu_2(\Delta)=2k+62$ $(k\geq2)$, then $d(S_1)=2$ with $R_{\lambda_1}(F)(y)=(y+1)^2\in\mathbb{F}_{\phi_1}[y]$. Let us replace $\phi_1$ by $x-(u+2^k)$. Then $F(x)=\cdots+A_2(x-(u+2^k))^2+A_1(x-(u+2^k))+A_0$, where $\nu_2(A_2)=6$,
	$$
	\begin{array}{lll}
		A_1&=&A+2^{k+1}a+2^k(72u^7)+2^{2k}(252u^6)+2^{3k}(504u^5)+2^{4k}(630u^4)+2^{5k}(504u^3)\\
		&&+2^{6k}(252u^2)+2^{7k}(72u)+2^{8k}(9)~~,\text{ and}\\
		A_0&=&B+2^kA+2^{2k}(a+36u^7)+2^{3k}(84u^6)+2^{4k}(126u^5)+2^{5k}(126u^4)+2^{6k}(84u^3)\\
		&&+2^{7k}(36u^2)+2^{8k}(9u)+2^{9k}~~.
	\end{array}
	$$ 
	Thus $\nu_2(A_1)=k+7$ and $\nu_2(A_0)\geq 2k+7$. We have 
	$$
	\begin{array}{lllll}
		A_0&\equiv&B+2^{2k}a\md{2^{2k+9}}&\equiv&\cfrac{2^{2k+6}}{3^{18}\cdot 7b_2^7}\left(\Delta_2+3^{18}\cdot 7a_2b_2^7\right)\md{2^{2k+9}}.
	\end{array}
	$$
	Since $\Delta_2+3^{18}\cdot 7a_2b_2^7\equiv \Delta_2-a_2b_2\md{8}$, then three cases arise:
	\item[-] If $\Delta_2\equiv a_2b_2+2\md{4}$, then $\nu_2(A_0)=2k+7$. Thus $N_{x-(u+2^k)}(F)=S_1+S_2$ has two sides joining $(0,2k+7)$, $(2,6)$, and $(9,0)$. It follows that $2\mathbb{Z}_K=\mathfrak{p}_1^2\mathfrak{p}_2^7$ with residue degree $1$ each ideal factor. Hence $\nu_2(i(K))=0$.
	
	\item[-] If $\Delta_2\equiv a_2b_2+4\md{8}$, then $\nu_2(A_0)=2k+8$ . Thus $N_{x-(u+2^k)}(F)=S_1+S_2$ has two sides joining $(0,2k+8)$, $(2,6)$, and $(9,0)$ with $d(S_2)=1$ and $R_{\lambda_1}(F)(y)=y^2+y+1$, which is irreducible over $\mathbb{F}_{x-(u+2^k)}$. It follows that $2\mathbb{Z}_K=\mathfrak{p}_1\mathfrak{p}_2^7$ with $f_1=2$ and $f_2=1$. Hence $\nu_2(i(K))=0$. 
	
	\item[-] If $\Delta_2\equiv a_2b_2\md{8}$, then $\nu_2(A_0)\geq2k+9$. Thus $N_{x-(u+2^k)}(F)=S_1+S_2+S_3$ has three sides joining $(0,\nu_2(A_0))$, $(1,k+7)$, $(2,6)$, and $(9,0)$. It follows that $2\mathbb{Z}_K=\mathfrak{p}_1\mathfrak{p}_2\mathfrak{p}_3^7$ with residue degree $1$ each ideal factor. Hence $2$ divides $i(K)$. Applying \cite[Corollary, p: 230]{En}, we get $\nu_2(i(K))=1$.	
	 
	\end{enumerate}
	\end{enumerate}
	
	\begin{flushright}
		$\square$
	\end{flushright}
	\textit{\textbf{Proof of Theorem \ref{thmp3}}}.\\ 
	If $\nu_3(ab)=0$, then since $\Delta=b(3^{18}b^7+2^{2}\cdot7^7a^9)$ is the discriminant of $F(x)$, thanks to the index formula $(1.1)$, $\nu_3((\mathbb{Z}_K:\mathbb{Z}[\alpha]))=0$ and so $\nu_3(i(K))=0$. Now assume that $3$ divides $ab$. Then we have the following cases:
	\begin{enumerate}
		\item If $a\equiv -1\md{3}$ and $\nu_3(b)\geq1$, then $F(x)\equiv x^2(x-1)(x^6+x^5+x^4+x^3+x^2+x+1) \md{3}$. Thus there are two prime ideals of $\mathbb{Z}_K$ lying above $3$ with residue degrees $1$ and $6$ provided by $x-1$ and $x^6+x^5+x^4+x^3+x^2+x+1$ respectively. On the other hand, $x$ can provides at most two prime ideals of $\mathbb{Z}_K$ lying above $3$ with residue degree $1$ each, or one prime ideal with residue degree $2$. Thus the possible factorizations of $3\mathbb{Z}_K$ are $\mathfrak{p}_1\mathfrak{p}_2\mathfrak{p}_3$ with $f_1=2$, $f_2=1$, and $f_3=6$, or $\mathfrak{p}_1^2\mathfrak{p}_2\mathfrak{p}_3$ with $f_1=f_2=1$ and $f_3=6$, or $\mathfrak{p}_{11}\mathfrak{p}_{12}\mathfrak{p}_{21}\mathfrak{p}_{31}$ with $f_{11}=f_{12}=f_{21}=1$, and $f_{31}=6$. In each case,  $\nu_3(i(K))=0$.
		\item If $a\equiv 1\md{3}$ and $\nu_3(b)\geq1$, then $F(x)\equiv x^2(x+1)(x^6-x^5+x^4-x^3+x^2-x+1)\md{3}$. Similar to the case above, $\nu_3(i(K))=0$.
		\item If $\nu_3(a)\geq 1$ and $b\equiv -1\md{3}$, then  $F(x)\equiv (x-1)^9 \md{3}$. Let $\phi=x-1$. Then $F(x)=\phi^9+9\phi^8+36\phi^7+84\phi^6+126\phi^5+126\phi^4+84\phi^3+(a+36)\phi^2+(2a+9)\phi+a+b+1$.
		\begin{enumerate}
			\item[(i)] If $(a,b)\in\{(0,2),(0,5),(3,-1),(3,2),(6,-1),(6,5)\}\md{9}$, then $N_{\phi}(F)=S_1$ has a single side of height $1$. Thus $3\mathbb{Z}_K=\mathfrak{p}_1^9$ with residue degree $1$. Hence $\nu_3(i(K))=0$.
			
			\item[(ii)] If $(a,b)\in\{(3,5),(6,2)\}\md{9}$, then $N_{\phi}(F)=S_1+S_2$ has two sides joining $(0,w)$, $(1,1)$, and $(9,0)$ with $w\geq2$. Thus $3\mathbb{Z}_K=\mathfrak{p}_1\mathfrak{p}_2^8$ with residue degree $1$ each ideal factor. Hence $\nu_3(i(K))=0$.
			
			\item[(iii)] If $(a,b)\in\{(0,8),(0,17),(9,-1),(9,8),(18,-1),(18,17)\}\md{27}$, then $N_{\phi}(F)=S_1+S_2$ has two sides joining $(0,2)$, $(3,1)$, and $(9,0)$. Thus $3\mathbb{Z}_K=\mathfrak{p}_1^3\mathfrak{p}_2^6$ with residue degree $1$ each ideal factor. Hence $\nu_3(i(K))=0$.
			
			\item[(iv)] If $(a,b)\in\{(0,-1),(18,8)\}\md{27}$, then $N_{\phi}(F)=S_1+S_2+S_3$ has three sides joining $(0,w)$, $(1,2)$, $(3,1)$, and $(9,0)$ with $w\geq3$. Thus the degree of each side $1$, and so $3\mathbb{Z}_K=\mathfrak{p}_1\mathfrak{p}_2^2\mathfrak{p}_3^6$ with residue degree $1$ each ideal factor. Hence $\nu_3(i(K))=0$.
			
			\item[(v)] If $(a,b)\in\{(9,17),(9,44),(36,17),(36,71),(63,44),(63,71)\}\md{81}$, then $N_{\phi}(F)=S_1+S_2$ has two sides joining $(0,3)$, $(3,1)$, and $(9,0)$. Thus $3\mathbb{Z}_K=\mathfrak{p}_1^3\mathfrak{p}_2^6$ with residue degree $1$ each ideal factor. Hence $\nu_3(i(K))=0$.
			
			\item[(vi)] If $(a,b)\equiv (9,71)\md{81}$, then we have the following sub-cases:
			\begin{enumerate}
				\item[(a)] If $a+b+1\equiv 81\md{243}$, then $N_{\phi}(F)=S_1+S_2$ has two sides joining $(0,4)$, $(3,1)$, and $(9,0)$ with $d(S_2)=1$ and  $R_{\lambda_1}(F)(y)=y^3-y^2+y+1$, which is irreducible over $\mathbb{F}_{\phi}$. It follows that $3\mathbb{Z}_K=\mathfrak{p}_1\mathfrak{p}_2^6$ with $f_1=3$ and $f_2=1$. Hence $\nu_3(i(K))=0$.
				
				\item[(b)] If $a+b+1\equiv 162\md{243}$, then $N_{\phi}(F)=S_1+S_2$ has two sides joining $(0,4)$, $(3,1)$, and $(9,0)$ with $d(S_2)=1$ and $R_{\lambda_1}(F)(y)=y^3-y^2+y-1=(y-1)(y^2+1)\in\mathbb{F}_{\phi}[y]$. Thus $3\mathbb{Z}_K=\mathfrak{p}_{11}\mathfrak{p}_{12}\mathfrak{p}_{21}^6$ with $f_{11}=f_{21}=1$ and $f_{12}=2$. Hence $\nu_3(i(K))=0$.
				
				\item[(c)] If $a+b+1\equiv 0\md{243}$, then $\nu_3(\Delta)\geq23$. Let $a_3=\cfrac{a}{9}$ and $u=\cfrac{2\cdot7^3a_3^4}{b^3}\in\mathbb{Q}$. Since $\nu_3(b^3)=0$, then $u\in\mathbb{Z}_3$. Let $\phi=x-u$. Then $F(x)=\cdots+84u^6\phi^3+(36u^7+a)\phi_2^2+A\phi_2+B$, with $\nu_3(36u^7+a)=2$,
				$$
				\begin{array}{lllllll}
					A&=&9u^8+2au&=&\cfrac{2^2\cdot 7^3a_3^5\Delta^3}{3^{52}b^{27}}-\cfrac{2^2\cdot 7^3a_3^5\Delta^2}{3^{33}b^{19}}+\cfrac{2^2\cdot 7^3a_3^5\Delta}{3^{15}b^{11}}~~, {\text{and}}\\
					&&&&\\
					B&=&u^9+au^2+b&=&\cfrac{2\Delta^4}{3^{71}\cdot 7b^{31}}-\cfrac{8\Delta^3}{3^{54}\cdot 7b^{23}}+\cfrac{4\Delta^2}{3^{35}\cdot 7b^{15}}+\cfrac{\Delta}{3^{18}\cdot 7b^{7}}~~.
				\end{array}
				$$ 
				Obviously, $\nu_3(A)=\nu_3(\Delta)-15$ and $\nu_3(B)=\nu_3(\Delta)-18$. It follows that $N_{\phi}(F)=S_1+S_2+S_3$ has three sides joining $(\nu_3(\Delta)-18)$, $(2,2)$, $(3,1)$, and $(9,0)$ with $d(S_2)=d(S_3)=1$. If $\nu_3(\Delta)$ is odd, then $d(S_1)=1$, and so $3\mathbb{Z}_K=\mathfrak{p}_1^2\mathfrak{p}_2\mathfrak{p}_3^6$ with residue degree $1$ each ideal factor. Hence $\nu_3(i(K))=0$.\\
				If $\nu_3(\Delta)$ is even, then $d(S_1)=2$ with $R_{\lambda_1}(F)(y)=-y^2+B_3\in\mathbb{F}_{\phi}[y]$ since $(36u^7+a)/9\equiv -1\md{3}$.
				\item[-] If $\Delta_3\equiv -b\md{3}$; $\Delta\equiv 1\md{3}$, then $R_{\lambda_1}(F)(y)=-y^2-1$, which is irreducible over $\mathbb{F}_{\phi}$. Thus $3\mathbb{Z}_K=\mathfrak{p}_1\mathfrak{p}_2\mathfrak{p}_3^6$ with $f_1=2$ and $f_2=f_3=1$. Hence $\nu_3(i(K))=0$.
				
				\item[-] If $\Delta_3\equiv -1\md{3}$, then $R_{\lambda_1}(F)(y)=-y^2+1=-(y-1)(y+1)\in\mathbb{F}_{\phi}[y]$. Thus $3\mathbb{Z}_K=\mathfrak{p}_{11}\mathfrak{p}_{12}\mathfrak{p}_{21}\mathfrak{p}_{31}^6$ with residue degree $1$ each ideal factor. Hence $3$ divides $i(K)$. Applying \cite[Corollary, p: 230]{En}, we get $\nu_3(i(K))=1$.
			\end{enumerate}
			
			\item[(vii)] If $(a,b)\equiv(36,44)\md{81}$, then we have the following cases:
			\begin{enumerate}
				\item[(a)] If $a+b+1\equiv 81\md{243}$, then $N_{\phi}(F)=S_1+S_2$ has two sides joining $(0,4)$, $(3,1)$, and $(9,0)$ with $d(S_2)=1$ and $R_{\lambda_1}(F)(y)=y^3-y^2+1$, which is irreducible over $\mathbb{F}_{\phi}$. It follows that $3\mathbb{Z}_K=\mathfrak{p}_1\mathfrak{p}_2^6$ with $f_1=3$ and $f_2=1$. Hence $\nu_3(i(K))=0$.
				
				\item[(b)] If $a+b+1\equiv 162\md{243}$, then $N_{\phi}(F)=S_1+S_2$ has two sides joining $(0,4)$, $(3,1)$, and $(9,0)$ with $d(S_2)=1$ and $R_{\lambda_1}(F)(y)=y^3-y^2-1=(y+1)(y^2+y-1)\in\mathbb{F}_{\phi}[y]$. It follows that $3\mathbb{Z}_K=\mathfrak{p}_{11}\mathfrak{p}_{12}\mathfrak{p}_{21}^6$ with $f_{11}=f_{21}=1$ and $f_{12}=2$. Hence $\nu_3(i(K))=0$.
				
				\item[(c)] If $a+b+1\equiv 0\md{243}$, then $\nu_3(\Delta)\geq23$. Let $a_3=\cfrac{a}{9}$ and $u=\cfrac{b}{20a_3^2}\in\mathbb{Q}$. Similar to the case (vi). (c) above, $3$ divides $i(K)$ if and only if $\nu_3(\Delta)$ is even and $\Delta_3\equiv -1\md{3}$. Also in this case, we have $\nu_3(i(K))=1$.
			\end{enumerate}
			
			\item[(viii)] If $(a,b)\equiv (63,17)\md{81}$, then we have the following sub-cases:
			\begin{enumerate}
				\item[(a)] If $a+b+1\equiv 81\md{243}$, then $\nu_3(\Delta)\geq23$. Let $a_3=\cfrac{a}{9}$ and $u=\cfrac{b}{20a_3^2}\in\mathbb{Q}$. Let $\phi=x-u$. Similar to the case (vi). (c) above, $3$ divides $i(K)$ if and only if $\nu_3(\Delta)$ is even and $\Delta_3\equiv -1\md{3}$. Also in this case, we have $\nu_3(i(K))=1$.
				
				\item[(b)] If $a+b+1\equiv 162\md{243}$, then $N_{\phi}(F)=S_1+S_2$ has two sides joining $(0,4)$, $(3,1)$, and $(9,0)$ with $d(S_2)=1$ and $R_{\lambda_1}(F)(y)=y^3-y^2-y-1$, which is irreducible over $\mathbb{F}_{\phi}$. Thus $3\mathbb{Z}_K=\mathfrak{p}_{1}\mathfrak{p}_{2}^6$ with $f_{1}=3$ and $f_{12}=2$. Hence $\nu_3(i(K))=0$.
				
				\item[(c)] If $a+b+1\equiv 0\md{243}$, then $N_{\phi}(F)=S_1+S_2+S_3$ has three sides joining $(0,w)$, $(1,3)$, $(3,1)$, and $(9,0)$ with $w\geq5$ and $R_{\lambda_2}(F)(y)=y^2-y-1$, which is irreducible over $\mathbb{F}_{\phi}$. Thus $3\mathbb{Z}_K=\mathfrak{p}_{1}\mathfrak{p}_2\mathfrak{p}_3^6$ with $f_1=f_3=1$ and $f_2=2$. Hence $\nu_3(i(K))=0$.
				
			\end{enumerate}
			
		\end{enumerate}
		\item If $\nu_3(a)\geq 1$ and $b\equiv 1\md{3}$, then  $F(x)\equiv (x+1)^9 \md{3}$. Let $\phi=x+1$. Then $F(x)=\phi^9-9\phi^8+36\phi^7-84\phi^6+126\phi^5-126\phi^4+84\phi^3+(a-36)\phi^2+(-2a+9)\phi+a+b-1$. The treatment of this case is similar to the case (3) above $(\nu_3(a)\geq 1$ and $b\equiv -1\md{3})$. The following are the obtained results:
		\begin{enumerate}
			\item[(i)] If $(a,b)\in\{(0,4),(0,7),(3,1),(3,4),(6,1),(6,7)\}\md{9}$, then $3\mathbb{Z}_K=\mathfrak{p}_1^9$ with residue degree $1$. Hence $\nu_3(i(K))=0$.
			
			\item[(ii)] If $(a,b)\in\{(3,7),(6,4)\}\md{9}$, then $3\mathbb{Z}_K=\mathfrak{p}_1\mathfrak{p}_2^8$ with residue degree $1$ each ideal factor. Hence $\nu_3(i(K))=0$.
			
			\item[(iii)] If $(a,b)\in\{(0,10),(0,19),(9,1),(9,10),(18,1),(18,19)\}\md{27}$, then $3\mathbb{Z}_K=\mathfrak{p}_1^3\mathfrak{p}_2^6$ with residue degree $1$ each ideal factor. Hence $\nu_3(i(K))=0$.
			
			\item[(iv)] If $(a,b)\in\{(0,1),(9,19)\}\md{27}$, then $3\mathbb{Z}_K=\mathfrak{p}_1\mathfrak{p}_2^2\mathfrak{p}_3^6$ with residue degree $1$ each ideal factor. Hence $\nu_3(i(K))=0$.
			
			\item[(v)] If $(a,b)\in\{(18,10),(18,37),(45,10),(45,64),(72,37),(72,64)\}\md{81}$, then $3\mathbb{Z}_K=\mathfrak{p}_1^3\mathfrak{p}_2^6$ with residue degree $1$ each ideal factor. Hence $\nu_3(i(K))=0$.
			
			\item[(vi)] If $(a,b)\equiv(18,64)\md{81}$, then we have the following sub-cases:
			\begin{enumerate}
				\item[(a)] If $a+b\equiv 82\md{243}$, then $3\mathbb{Z}_K=\mathfrak{p}_1\mathfrak{p}_2^6$ with $f_1=3$ and $f_2=1$. Hence $\nu_3(i(K))=0$.
				
				\item[(b)] If $a+b\equiv 1\md{243}$, then $3\mathbb{Z}_K=\mathfrak{p}_1\mathfrak{p}_2\mathfrak{p}_3^6$ with $f_1=f_3=1$ and $f_2=2$. Hence $\nu_3(i(K))=0$.
				
				\item[(c)] If $a+b\equiv 163\md{243}$ and $\nu_3(\Delta)$ is odd, then $3\mathbb{Z}_K=\mathfrak{p}_1^2\mathfrak{p}_2\mathfrak{p}_3^6$ with residue degree $1$ each ideal factor. Hence $\nu_3(i(K))=0$.
				
				\item[(d)] If $a+b\equiv 163\md{243}$, $\nu_3(\Delta)$ is even, and $\Delta_3\equiv 1\md{3}$, then $3\mathbb{Z}_K=\mathfrak{p}_1\mathfrak{p}_2\mathfrak{p}_3^6$ with $f_1=2$ and $f_2=f_3=1$. Hence $\nu_3(i(K))=0$.
				
				\item[(e)] If $a+b\equiv 163\md{243}$, $\nu_3(\Delta)$ is even, and $\Delta_3\equiv -1\md{3}$, then $3\mathbb{Z}_K=\mathfrak{p}_1\mathfrak{p}_2\mathfrak{p}_3\mathfrak{p}_3^6$ with residue degree $1$ each ideal factor. Hence $\nu_3(i(K))=1$.
			\end{enumerate}
		
			\item[(vii)] If $(a,b)\in\{(45,37),(72,10)\}\md{81}$, then we have the following sub-cases:
			
			\begin{enumerate}
				\item[(a)] If $a+b\equiv 82\md{243}$, then $3\mathbb{Z}_K=\mathfrak{p}_1\mathfrak{p}_2\mathfrak{p}_3^6$ with $f_1=f_3=1$ and  $f_2=2$. Hence $\nu_3(i(K))=0$.
				
				\item[(b)] If $a+b\equiv 163\md{243}$, then $3\mathbb{Z}_K=\mathfrak{p}_1\mathfrak{p}_3^6$ with $f_1=3$ and $f_2=1$. Hence $\nu_3(i(K))=0$.
				
				\item[(c)] If $a+b\equiv 1\md{243}$ and $\nu_3(\Delta)$ is odd, then $3\mathbb{Z}_K=\mathfrak{p}_1^2\mathfrak{p}_2\mathfrak{p}_3^6$ with residue degree $1$ each ideal factor. Hence $\nu_3(i(K))=0$.
				
				\item[(d)] If $a+b\equiv 1\md{243}$, $\nu_3(\Delta)$ is even, and $\Delta_3\equiv 1\md{3}$, then $3\mathbb{Z}_K=\mathfrak{p}_1\mathfrak{p}_2\mathfrak{p}_3^6$ with $f_1=2$ and $f_2=f_3=1$. Hence $\nu_3(i(K))=0$.
				
				\item[(e)] If $a+b\equiv 1\md{243}$, $\nu_3(\Delta)$ is even, and $\Delta_3\equiv -1\md{3}$, then $3\mathbb{Z}_K=\mathfrak{p}_1\mathfrak{p}_2\mathfrak{p}_3\mathfrak{p}_3^6$ with residue degree $1$ each ideal factor. Hence $\nu_3(i(K))=1$.
			\end{enumerate}
		\end{enumerate}		
		
		\item If $\nu_3(a)\geq 1$ and $\nu_3(b)\geq1$, then  $F(x)\equiv x^9 \md{3}$. Let $\phi=x$. By assumption, $\nu_3(a)\leq6$ or $\nu_3(b)\leq8$.
		\begin{enumerate}
			\item[(i)] If $\nu_3(b)<\cfrac{9}{7}\nu_3(a)$, then $N_{\phi}(F)=S_1$ has a single side joining $(0,\nu_3(b))$ and $(9,0)$ with $d(S_1)\in\{1,3\}$ since $\nu_3(b)\leq 8$ by assumption. Thus there are at most three prime ideals of $\mathbb{Z}_K$ lying above $3$ with residue degree $1$ each ideal factor. Hence $\nu_3(i(K))=0$.
			
			\item[(ii)] If $\nu_3(b)>\cfrac{9}{7}\nu_3(a)$, then $N_{\phi}(F)=S_1+S_2$ has two sides joining $(0,\nu_3(b))$, $(2,\nu_3(a))$, and $(9,0)$ with $d(S_1)\in\{1,2\}$ and $d(S_2)=1$ since $\nu_3(a)\leq 6$ by assumption. Thus there are at most three prime ideals of $\mathbb{Z}_K$ lying above $3$ with residue degree $1$ each ideal factor. Hence $\nu_3(i(K))=0$.
		\end{enumerate}
		
	\end{enumerate}
	
	\begin{flushright}
		$\square$
	\end{flushright}
	
	\textit{\textbf{Proof of Theorem \ref{thmp}}}.\\
	For $p=5$. Since $\Delta=b(3^{18}b^7+2^{2}\cdot7^7a^9)$, then $\nu_5(\Delta)\geq1$ if and only if $\nu_5(b)\geq1$ or $(a,b)\in\{(1,3),(2,4),(3,1),(4,2)\}\md{5}$.
	\begin{enumerate}
		\item If $(a,b)\in\{(1,3),(2,4),(3,1),(4,2)\}\md{5}$ or $\nu_5(b)\geq1$ and $\nu_5(a)=0$, then one can check easily that $F(x)\equiv \phi_{i}^2\cdot\phi_{j}\cdot\phi_{k}\md{5}$ with $\deg(\phi_{i})=\deg(\phi_{j})=1$, $\deg(\phi_{k})=6$, and the polynomials $\phi_{i}$, $\phi_{j}$, and $\phi_{k}$ are coprime irreducible over $\mathbb{F}_5$ for every $i,j,k$. Thus there are at most three prime ideals of $\mathbb{Z}_K$ lying above $5$ with residue degree $1$ each ideal factor. Hence $\nu_5(i(K))=0$.
		
		\item If $\nu_5(b)\geq 1$, and $\nu_5(a)\geq1$ then  $F(x)\equiv x^9 \md{5}$. Let $\phi=x$.  By assumption, $\nu_5(a)\leq6$ or $\nu_5(b)\leq8$.
		\begin{enumerate}
			\item[(i)] If $\nu_5(b)<\cfrac{9}{7}\nu_5(a)$, then $N_{\phi}(F)=S_1$ has a single side joining $(0,\nu_5(b))$ and $(9,0)$ with $d(S_1)\in\{1,3\}$ since $\nu_5(b)\leq 8$ by assumption. Thus there are at most three prime ideals of $\mathbb{Z}_K$ lying above $5$ with residue degree $1$ each ideal factor. Hence $\nu_5(i(K))=0$.
			
			\item[(ii)] If $\nu_5(b)>\cfrac{9}{7}\nu_5(a)$, then $N_{\phi}(F)=S_1+S_2$ has two sides joining $(0,\nu_5(b))$, $(5,\nu_5(a))$, and $(9,0)$ with $d(S_1)\in\{1,2\}$ and $d(S_2)=1$ since $\nu_5(a)\leq 6$ by assumption. Thus there are at most three prime ideals of $\mathbb{Z}_K$ lying above $5$ with residue degree $1$ each ideal factor. Hence $\nu_5(i(K))=0$.
		\end{enumerate}
	\end{enumerate}
	For $p=7$. Since $\Delta=b(3^{18}b^7+2^{2}\cdot7^7a^9)$, thanks to the index formula, if $\nu_7(b)=0$, then $\nu_7((\mathbb{Z}_K:\mathbb{Z}[\alpha]))=0$, and so $\nu_7(i(K))=0$. Assume that $7$ divides $b$. Then we have the following cases:
	\begin{enumerate}
		\item If $\nu_7(a)=0$, then $F(x)\equiv x^2(x+a)^7\md{7}$. Let $\phi_1=x$ and $\phi_2=x+a$. Then
		$$
		\begin{array}{lll}
			F(x)&=&\phi_1^9+a\phi_1^2+b,\\
			&=&\cdots+36a^2\phi_2^7-84a^3\phi_2^6+126a^4\phi_2^5-126a^5\phi_2^4+84a^6\phi_2^3+(-36a^7+a)\phi_2^2\\
			&&+(-2a^2+9a^8)\phi_2+b+a^3-a^9.
		\end{array}
		$$ 
		Obviously, $\phi_1$ can provides at most two prime ideals of $\mathbb{Z}_K$ lying above $7$ with residue degree $1$ each ideal factor.
		On the other hand, $\nu_7(84a^6)=1$, then $\phi_2$ provides a prime ideal $\mathfrak{p}$ of $\mathbb{Z}_K$ lying above $7$ with ramification index $e(\mathfrak{p}/7)\geq4$. Thus $\phi_2$ can provides at most four prime ideals of $\mathbb{Z}_K$ lying above $7$ with residue degree $1$ each ideal factor. We conclude that there are at most six prime ideals of $\mathbb{Z}_K$ lying above $7$ with residue degree $1$ each ideal factor. Hence $\nu_7(i(K))=0$.
		
		\item If $\nu_7(a)\geq1$, then $F(x)\equiv x^{9}\md{7}$. Let $\phi=x$. By assumption, $\nu_7(a)\leq6$ or $\nu_7(b)\leq8$.
		\begin{enumerate}
			\item[(i)] If $\nu_7(b)<\cfrac{9}{7}\nu_7(a)$, then $N_{\phi}(F)=S_1$ has a single side joining $(0,\nu_7(b))$ and $(9,0)$ with $d(S_1)\in\{1,3\}$ since $\nu_7(b)\leq 8$ by assumption. Thus there are at most three prime ideals of $\mathbb{Z}_K$ lying above $7$ with residue degree $1$ each ideal factor. Hence $\nu_7(i(K))=0$.
			
			\item[(ii)] If $\nu_7(b)>\cfrac{9}{7}\nu_7(a)$, then $N_{\phi}(F)=S_1+S_2$ has two sides joining $(0,\nu_7(b))$, $(2,\nu_7(a))$, and $(9,0)$ with $d(S_1)\in\{1,2\}$ and $d(S_2)=1$ since $\nu_7(a)\leq 6$ by assumption. Thus there are at most three prime ideals of $\mathbb{Z}_K$ lying above $7$ with residue degree $1$ each ideal factor. Hence $\nu_7(i(K))=0$.
		\end{enumerate}
		
	\end{enumerate}
	
	For $p\geq11$, since there is at most $9$ prime ideals of $\mathbb{Z}_{K}$ lying above $p$ with residue degree $1$ each ideal factor, and there is at least $p\geq 11$ monic irreducible polynomial of degree $f$ in $\mathbb{F}_p[x]$ for every positive integer $f$, we conclude that $p$ does not divide $i(K)$.
	\begin{flushright}
		$\square$
	\end{flushright}

	\section{Examples}
	Let $F(x)=x^9+ax^2+b\in \mathbb{Z}[x]$ be a monic irreducible polynomial and $K$ the nonic number field generated by a complex root of $F(x)$.
	\begin{enumerate}
		\item For $a=54$ and $b=87$, we have $\nu_2(b)=0$, then by Theorem $\ref{thmp2}$, $\nu_2(i(K))=0$. On the other hand, $(a,b)\equiv(54,6)\md{81}$, then by Theorem $\ref{thmp3}$, $\nu_3(i(K))=0$. We conclude that $i(K)=1$. 
		
		\item For $a=64$ and $b=256$, we have $a\equiv 64\md{128}$, $b\equiv 256\md{256}$, $\nu_2(\Delta)=70$, and $\Delta_2\equiv a_2b_2\equiv 1\md{8}$, then by Theorem $\ref{thmp2}$, $\nu_2(i(K))=1$. On the other hand, $(a,b)\equiv (64,13)\md{81}$, then by Theorem $\ref{thmp3}$, $\nu_3(i(K))=0$. Hence $i(K)=2$.

		\item For $a=90$ and $b=19835$, we have $\nu_2(b)=0$, by Theorem $\ref{thmp2}$, $\nu_2(i(K))=0$. On the other hand, $(a,b)\equiv (9,71)\md{81}$, $a+b\equiv -\md{243}$, $\nu_3(\Delta)=26$, and $\Delta_3\equiv -1\md{3}$, by Theorem $\ref{thmp3}$, $\nu_3(i(K))=1$. Hence $i(K)=3$.

		\item For $a=99$ and $b=8055028$, we have $\nu_2(a)=0$, $\nu_2(b)=2$, and $b_2\equiv -a\equiv 5\md{8}$, then by Theorem $\ref{thmp2}$, $\nu_2(i(K))=1$. On the other hand, $(a,b)\equiv(9,50)\md{81}$, $a+b\equiv 163\md{243}$, $\nu_3(\Delta)=28$, and $\Delta_3\equiv -1\md{3}$, then by Theorem $\ref{thmp3}$, $\nu_3(i(K))=1$. Hence $i(K)=6$.
		
	\end{enumerate}
	\text{ }\newline
	{\bf Conflict of interest}\\
	Not Applicable.\\
	{\bf Data availability}\\
	Not applicable.\\
	{\bf Author Contribution and Funding Statement}\\
	Not applicable.

\end{document}